
\documentstyle{amsppt}
\baselineskip18pt
\magnification=\magstep1
\pagewidth{30pc}
\pageheight{45pc}
\hyphenation{co-deter-min-ant co-deter-min-ants pa-ra-met-rised
pre-print pro-pa-gat-ing pro-pa-gate
fel-low-ship Cox-et-er dis-trib-ut-ive}
\def\leaderfill{\leaders\hbox to 1em{\hss.\hss}\hfill}

\

\def\idest{i.e.,\ }

\def\a{{\alpha}}
\def\be{{\beta}}
\def\g{{\gamma}}

\def\d{{\delta}}

\def\e{{\varepsilon}}

\def\br{{\bold r}}
\def\bu{{\bold u}}
\def\bv{{\bold v}}

\def\bx{{\bold x}}

\def\b0{\text{\bf 0}}

\def\boxit#1{\vbox{\hrule\hbox{\vrule \kern3pt
\vbox{\kern3pt\hbox{#1}\kern3pt}\kern3pt\vrule}\hrule}}
\def\rabbit{\vbox{\hbox{\kern0pt
\vbox{\kern0pt{\hbox{---}}\kern3.5pt}}}}

\def\tableau#1{
        \hbox {
                \hskip -10pt plus0pt minus0pt
                \raise\baselineskip\hbox{
                \offinterlineskip
                \hbox{#1}}
                \hskip0.25em
        }
}

\def\tabCol#1{
\hbox{\vtop{\hrule
\halign{\strut\vrule\hskip0.5em##\hskip0.5em\hfill\vrule\cr\lower0pt
\hbox\bgroup$#1$\egroup \cr}
\hrule
} } \hskip -10.5pt plus0pt minus0pt}

\def\CR{
        $\egroup\cr
        \noalign{\hrule}
        \lower0pt\hbox\bgroup$
}



\def\blank#1#2{
\hbox to #1{\hfill \vbox to #2{\vfill}}
}


\def\strut{\vrule height10pt depth5pt width0pt}

\def\bs{{\bold s}}

\topmatter
\title Freely braided elements in Coxeter groups
\endtitle

\author R.M. Green and J. Losonczy \endauthor
\affil Department of Mathematics and Statistics\\ Lancaster
University\\ Lancaster LA1 4YF\\ England\\ {\it  E-mail:}
r.m.green\@lancaster.ac.uk\\
\newline
Department of Mathematics\\ Long Island University\\ Brookville,
NY  11548\\ USA\\ {\it  E-mail:} losonczy\@e-math.ams.org\\
\newline
\endaffil

\abstract We introduce a notion of ``freely braided element" for
simply laced Coxeter groups. We show that an arbitrary group
element $w$ has at most $2^{N(w)}$ commutation classes of reduced
expressions, where $N(w)$ is a certain statistic defined in terms
of the positive roots made negative by $w$.  This bound is
achieved if $w$ is freely braided. In the type $A$ setting, we
show that the bound is achieved only for freely braided $w$.
\endabstract



\endtopmatter

\head 1. Introduction \endhead

A well-known result in the theory of Coxeter groups states that
any two reduced expressions for the same element $w$ of a Coxeter
group are equivalent under the equivalence relation generated by
braid relations.  If $w$ is such that any two of its reduced
expressions are equivalent by short braid relations (\idest
iterated commutations of commuting generators), we call $w$
``fully commutative," following Stembridge \cite{8}.

In this paper, we introduce and study ``freely braided elements''
for an arbitrary simply laced Coxeter group.  This is a more
general class of elements than the aforementioned fully
commutative elements.  The idea behind the definition is that
although it may be necessary to use long braid relations in order
to pass between two reduced expressions for a freely braided
element, the necessary long braid relations in a certain sense do
not interfere with one another.

Every reduced expression for a Coxeter group element $w$
determines a total ordering of the set of positive roots made
negative by $w$. The resulting sequences, called root sequences,
play a central role in this paper.  We are particularly interested
in triples of roots of the form $\{\a, \a + \be, \be\}$, where
each root in the triple is made negative by $w$ (inversion
triples). An inversion triple that occurs consecutively in some
root sequence for $w$ will be called contractible, and the
statistic $N(w)$ mentioned in the abstract above is the number of
contractible inversion triples of $w$. We note that Fan and
Stembridge have already shown in \cite{5, Theorem 2.4} that an
element of a simply laced Coxeter group is fully commutative if
and only if it has no inversion triples.

In \S4 of this paper, we prove that the number of commutation
classes (short braid equivalence classes of reduced expressions)
of $w$ is bounded above by $2^{N(w)}$, and this bound is achieved
if $w$ is freely braided. Furthermore, a freely braided element
$w$ has a root sequence in which each contractible inversion
triple occurs as a consecutive subword (Theorem 4.2.3). In \S5, we
prove that if the Coxeter group is of type $A$, then the
$2^{N(w)}$ bound is achieved only for freely braided $w$. Our
proof relies on the fact that every inversion triple is
contractible in the type $A$ setting.

Along the way, we give a short proof that the commutation graph of
any element of a simply laced Coxeter group is bipartite,
extending a result in \cite{4}.

Another possible approach to proving these results is to use
Viennot's heaps of pieces \cite{10}.  These are certain labelled
posets that can be associated to commutation classes of reduced
expressions for elements of a Coxeter group (see \cite{8, \S1.2}).
It can be shown that the dual of the heap of a reduced expression
is isomorphic, as an abstract poset, to a certain poset arising
naturally from the associated root sequence. However, we find it
more convenient to argue directly with root sequences in this
paper.


\head 2. Preliminaries \endhead

\subhead 2.1 Basic terminology and notation \endsubhead

Let $(W,S)$ be a Coxeter system, with Coxeter matrix $\{ m(s,t)
\}_{s,t\in S}$. It will be assumed throughout this paper that
$(W,S)$ is simply laced, so that $m(s,t) \in \{2,3\}$ for all
pairs of distinct $s,t\in S$. We refer to $W$ itself as a
``Coxeter group". The basic facts concerning Coxeter systems can
be found in \cite{3} or \cite{6}.

Let $S^*$ be the free monoid generated by $S$.  The elements of
$S^*$ will be written as finite sequences, e.g., $\bs =
(s_1,\dots,s_n)$. We define the {\it length} of any $\bs \in S^*$
to be the number of its entries.

There is a natural morphism of monoids $\phi: S^* \longrightarrow
W$, given by $\phi(s_1,\dots,s_n) = s_1\cdots s_n$ (the empty
sequence is mapped to the identity, $e$). We say that an element
$\bs \in S^*$ {\it represents} its image $w=\phi(\bs)\in W$;
furthermore, if the length of $\bs$ is minimal among the lengths
of all the sequences that represent $w$, then we say that $\bs$ is
{\it reduced}, and we call $\bs$ a {\it reduced expression} for
$w$. The {\it length} of $w$, denoted by $\ell(w)$, is then equal
to the length of $\bs$.

Let $V$ be a vector space over ${\Bbb R}$ with basis $\{\a_s :
s\in S\}$, and let $B$ be the Coxeter form on $V$ associated to
$(W,S)$. This is the symmetric bilinear form satisfying
$B(\a_s,\a_t) = -\cos{\frac{\pi}{m(s,t)}}$ for all $s,t\in S$. We
shall view $V$ as the underlying space of a reflection
representation of $W$, determined by the equalities $s\,\a_t =
\a_t - 2B(\a_t,\a_s)\a_s$, for $s,t\in S$. The Coxeter form is
preserved by $W$ relative to this representation.

Define $\Phi = \{w\a_s : w\in W\text{ and }s\in S\}$.  This is the
root system of $(W,S)$.  Let $\Phi^+$ be the set of all $\be \in
\Phi$ such that $\be$ is expressible as a linear combination of
the $\a_s$ with nonnegative coefficients, and let $\Phi^- =
-\Phi^+$. Then $\Phi$ equals the disjoint union $\Phi^+ \cup
\Phi^-$ \cite{6, \S5.4}. The elements of $\Phi^+$ (respectively,
$\Phi^-$) are called ``positive" (respectively, ``negative")
roots. The $\a_s$ are often referred to as ``simple" roots.

One can associate to each $w\in W$ the set $\Phi(w) = \Phi^+ \cap
w^{-1}(\Phi^-)$, which we call the {\it inversion set} of $w$. It
is well known that $\Phi(w)$ has $\ell(w)$ elements and that $w$
is uniquely determined by $\Phi(w)$. Given any reduced expression
$(s_1,\dots,s_n)$ for $w$, one has $\Phi(w) = \{r_1,\dots,r_n \}$,
where $r_1 = \a_{s_n}$ and $r_i = s_n\cdots
s_{n-i+2}(\a_{s_{n-i+1}})$ for all $i\in \{2,\dots,n\}$ \cite{6,
Exercise 1, \S5.6}.  Following \cite{5}, we form the sequence $\br
=(r_1,\dots,r_n)$, and call it the {\it root sequence of }
$(s_1,\dots,s_n)$, or a {\it root sequence for } $w$. Notice that
any initial segment of a root sequence is also a root sequence for
some element of $W$.

Let $w\in W$. Suppose that $\a,\be$ are positive roots and $a,b$
are positive integers such that $a\a+b\be$ is a positive root.
Then one has the following two implications: $\a, \be \in \Phi(w)
\Rightarrow a\a+b\be \in \Phi(w)$ and $\a, \be \notin \Phi(w)
\Rightarrow a\a+b\be \notin \Phi(w)$. This property of inversion
sets is sometimes referred to as ``biconvexity" (cf. \cite{2,
\S3}).

\subhead 2.2 Inversion triples \endsubhead

Assume that $(W,S)$ is of arbitrary simply laced type, unless
stated otherwise.

\definition{Definition 2.2.1}
Let $w \in W$.  Any subset of $\Phi(w)$ of the form $\{ \a, \a +
\be, \be \}$ will be called an {\it inversion triple} of $w$. We
say that an inversion triple $T$ of $w$ is {\it contractible} if
there is a root sequence for $w$ in which the elements of $T$
appear consecutively (in some order).  If the contractible
inversion triples of $w$ are pairwise disjoint, then $w$ is said
to be {\it freely braided}.
\enddefinition

\remark{Remark 2.2.2} Let $\{ \a, \a + \be, \be \}$ be an
inversion triple of $w \in W$. It follows immediately from the
final two paragraphs of \S2.1 that $\a+\be$ must occur between
$\a$ and $\be$ in any root sequence for $w$.
\endremark

\vskip 5pt

The following example establishes the existence of
non-contractible inversion triples.

\example{Example 2.2.3} Suppose that $(W,S)$ has a parabolic
subgroup of type $D_4$, generated by $s_1,s_2,s_3,s_4\in S$, where
$s_2s_i$ has order $3$ for all $i\neq 2$, and $s_is_j$ has order
$2$ whenever $i,j \in \{1,3,4\}$ are distinct. Write $\a_i$ for
the simple root $\a_{s_i}$. Consider the element $w =
s_2s_1s_3s_4s_2s_4s_3s_1s_2$. The root sequence corresponding to
the reduced expression $(s_2,s_1,s_3,s_4,s_2,s_4,s_3,s_1,s_2)$ is
$$\eqalign{ (\a_2,\ \a_1+\a_2,\ &\a_2+\a_3,\ \a_2+\a_4,\
\a_1+2\a_2+\a_3+\a_4,\ \cr &\a_1+\a_2+\a_3,\  \a_1+\a_2+\a_4,\
\a_2+\a_3+\a_4,\ \a_1+\a_2+\a_3+\a_4) .\cr}$$  Observe that $\a_2$
belongs to precisely one inversion triple of $w$, namely $$
\{\a_2,\ \a_1+2\a_2+\a_3+\a_4,\ \a_1+\a_2+\a_3+\a_4\} .$$ However,
there is no root sequence for $w$ in which these roots appear
consecutively.  This can be seen as follows.  By inspection of the
above root sequence, $\a_2$ is the only simple root in $\Phi(w)$.
Hence, $\a_2$ must appear first in any root sequence for $w$. The
second root in any root sequence for $w$ must be a linear
combination of at most two simple roots; this rules out
$\a_1+2\a_2+\a_3+\a_4$ as a possibility.
\endexample

It will be shown in \S5.1 that in the type $A$ setting, every
inversion triple is contractible.

\head 3. Braid relations  \endhead

In this section, $(W,S)$ is assumed to be of arbitrary simply
laced type.  As in \S2, we denote the associated Coxeter matrix by
$\{m(s,t)\}_{s,t\in S}$.

\subhead 3.1 More on root sequences
\endsubhead

The group $W$ is by definition generated by the elements of $S$,
subject only to the relations $(st)^{m(s,t)} = e$, for $s,t \in S$
\cite{6, \S5}. These relations can be restated as $s^2 = e$ for
$s\in S$, $st = ts$ if $m(s,t)=2$, and $sts=tst$ if $m(s,t)=3$.
Given any $s,t \in S$ and any nonnegative integer $m$, we write
$(s,t)_m$ for the length $m$ sequence $(s,t,s,\dots)\in S^*$. Let
$\approx$ be the equivalence relation on $S^*$ generated by the
elementary relations $(s,t)_{m(s,t)}\approx (t,s)_{m(s,t)}$, for
$s,t\in S$. When $m(s,t)\neq 1$, we call such an elementary
relation a {\it braid relation}, qualifying it {\it short} or {\it
long}, according as $m(s,t) = 2$ or $3$.

Let $w\in W$.  A well-known result, attributed variously to
Matsumoto \cite{7} and Tits \cite{9}, states that the equivalence
class relative to $\approx$ of any reduced expression for $w$
coincides with the set of all reduced expressions for $w$.

We use the term {\it commutation class} to refer to any
equivalence class arising from the equivalence relation on $S^*$
generated by the short braid relations.

Applying a braid relation to a given reduced expression
corresponds to applying a permutation to the root sequence
associated with that reduced expression. The following proposition
makes this correspondence explicit.

\proclaim{Proposition 3.1.1} Let $w \in W$, let $\bu, \bv \in S^*$
and let $s,s',s''\in S$. Denote the length of $\bv$ by $k$.
\item{\rm (a)}{Assume that $(\bu,s,s',\bv)$ is a reduced expression
for $w$, and let $\br = (r_i)$ be the root sequence of
$(\bu,s,s',\bv)$.
\item{\rm (i)}{Suppose $m(s,s') = 2$, so that
$(\bu,s',s,\bv)$ is also a reduced expression for $w$. Then the
root sequence $\br'$ of $(\bu,s',s,\bv)$ can be obtained from
$\br$ by interchanging $r_{k+1}$ and $r_{k+2}$, which are mutually
orthogonal.}
\item{\rm (ii)}{If $r_{k+1}$ and $r_{k+2}$ are orthogonal,
then $m(s,s') = 2$.} }
\item{\rm (b)}{Assume that $(\bu,s,s',s'',\bv)$ is a reduced
expression for $w$, and let $\br = (r_i)$ be the root sequence of
$(\bu,s,s',s'',\bv)$.
\item{\rm (i)}{Suppose $s'' = s$, so that $m(s,s')=3$ and
$(\bu,s',s,s',\bv)$ is also a reduced expression for $w$. Then the
root sequence $\br'$ of $(\bu,s',s,s',\bv)$ can be obtained from
$\br$ by interchanging $r_{k+1}$ and $r_{k+3}$. Furthermore, we
have $r_{k+1} + r_{k+3}=r_{k+2}$.}
\item{\rm (ii)}{If $r_{k+1} + r_{k+3} = r_{k+2}$, then
$s'' = s \ne s'$ and $m(s,s') = 3$.} }
\endproclaim

\demo{Proof} We first prove the proposition under the additional
hypothesis that $\bv$ has length $0$. For simplicity of notation,
we denote the simple roots corresponding to $s$, $s'$ and $s''$ by
$\a$, $\a'$ and $\a''$, respectively.

Concerning (a), if $m(s,s')=2$, then the reduced expression
$(\bu,s,s')$ has root sequence $(\a',\a,r_3,\dots)$ and
$(\bu,s',s)$ has root sequence $(\a,\a',r_3,\dots)$. Statement (i)
immediately follows. For (ii), we observe that if the roots $\a'$
and $s'(\a)$ are orthogonal, then so are $\a$ and $\a'$, and hence
$m(s,s')=2$.

We turn to (b). If $s=s''$, then the expression $(\bu,s,s',s)$ is
reduced and has root sequence $(\a,\a +\a',\a',r_4,\dots)$;
likewise, $(\bu,s',s,s')$ is reduced and has root sequence
$(\a',\a + \a',\a,r_4,\dots)$. This establishes (i). To prove
(ii), we observe that in order for the equation
$$ \a'' + s''s'(\a) = s''(\a')
$$ to hold, we must have $s'' = s$; otherwise, the simple root
$\a$ would lie in the support of the left hand side but not the
right hand side (since $s\neq s'$). It follows that $m(s,s')=3$ if
$r_1+r_3=r_2$.

The general case (where $\bv$ has length $\geq 0$) follows from
the above, using the definition of root sequence together with the
fact that $W$ acts linearly on $V$ and preserves $B$. \qed\enddemo

\definition{Definition 3.1.2} Let $\br$ and $\br'$ be as
in part (a)(i) (respectively, part (b)(i)) of Proposition 3.1.1.
We call the passage from $\br$ to $\br'$ a {\it short braid move}
(respectively, {\it long braid move}). We say that two root
sequences are {\it commutation equivalent} if one can be
transformed into the other by applying a (possibly empty) sequence
of short braid moves.  The set of all root sequences that are
commutation equivalent to a given root sequence will be called the
{\it commutation class} of that root sequence.\enddefinition

\remark{Remark 3.1.3} Let $w\in W$. The recipe that associates a
root sequence to a reduced expression defines a bijection from the
set of all reduced expressions for $w$ to the set of all root
sequences for $w$. By Proposition 3.1.1, this bijection is
compatible with the application of both long and short braid
moves. Hence, by the result of Matsumoto and Tits cited earlier,
any root sequence for $w$ may be transformed into any other by
applying a (possibly empty) sequence of long and short braid
moves.
\endremark

\definition{Definition 3.1.4}
Let $w \in W$ and let $\br = (r_1,\dots, r_n)$ be any root
sequence for $w$.  Following \cite{4}, we introduce a partial
ordering $\leq_\br$ of $\Phi(w)$ by stipulating that $r_i <_\br
r_j$ if $i < j$ and the roots $r_i$ and $r_j$ are not orthogonal
relative to $B$.
\enddefinition

The following result is closely related to \cite{5, Proposition
2.2}.

\proclaim{Proposition 3.1.5} Let $w \in W$ and let $\br$ and
$\br'$ be root sequences for $w$ with respective partial orderings
$\leq_\br$ and $\leq_{\br'}$ of $\Phi(w)$. Then $\leq_\br$ and
$\leq_{\br'}$ agree if and only if $\br$ and $\br'$ are
commutation equivalent.
\endproclaim

\demo{Proof} For the ``if\," part, it is enough to treat the case
where $\br$ and $\br'$ differ by a single short braid move. In
such a situation, the conclusion follows immediately from
Definition 3.1.4.

For the converse, write $\br = (r_1,\dots,r_n)$ and $\br' =
(r'_1,\dots,r'_n)$.  Let $g$ be the permutation of $\{1,\dots,n\}$
such that $r_i = r'_{g(i)}$ for all $i \in \{1,\dots,n\}$. If
$\leq_\br$ and $\leq_{\br'}$ agree, then, for all $i,j \in
\{1,\dots,n\}$ such that $r'_{g(i)}$ and $r'_{g(j)}$ are
nonorthogonal, we have $i < j$ if and only if $g(i) < g(j)$. In
particular, this means that if $i$ is an index less than $n$ such
that $r'_{g(i)}$ is not orthogonal to $r'_{g(n)}$, then $g(i) <
g(n)$. Hence, $r'_{g(n)}$ can be moved to the last entry of the
root sequence $\br'$ using only short braid moves. Denote the
resulting root sequence by $\br''$.  Since any initial segment of
a root sequence is also a root sequence, and since any element of
$W$ is uniquely determined by its inversion set, there is a $w'
\in W$ such that the first $n-1$ entries of $\br$ and the first
$n-1$ entries of $\br''$ are both root sequences for $w'$. The
partial orders corresponding to these shorter root sequences
agree. By induction, each of the shorter root sequences can be
transformed into the other using only short braid moves.
\qed\enddemo

\subhead 3.2 More on inversion triples \endsubhead

The proposition below characterizes contractible inversion
triples. Its proof depends on the following basic fact from the
theory of Coxeter groups: given $w\in W$ and $s\in S$, one has
$\a_s \in \Phi(w)$ if and only if $w$ has a reduced expression
with $s$ as its last entry (see \cite{6, Proposition 5.7, Theorem
5.8}).

\proclaim{Proposition 3.2.1} Let $w \in W$ and let $T = \{\a, \a +
\be, \be\}$ be an inversion triple of $w$. The following are
equivalent\,{\rm :}
\itemitem{\rm (i)}{$T$ is contractible\,{\rm ;}}
\itemitem{\rm (ii)}{there is a root sequence $\br$ for $w$
such that $\a+\be$ covers $\a$ or $\be$ relative to $\leq_\br${\rm
;}}
\itemitem{\rm (iii)}{there is a root sequence $\br$ for $w$
such that $\a$ or $\be$ covers $\a+\be$ relative to $\leq_\br$.}
\endproclaim

\demo{Proof} The implications (i) $\Rightarrow$ (ii) and (i)
$\Rightarrow$ (iii) are immediate from the definition of
contractible inversion triple.

To prove (ii) $\Rightarrow$ (i), it is enough by symmetry to treat
the case where there is a root sequence $\br$ for $w$ such that
$\a + \be$ covers $\a$ relative to $\leq_\br$.  We may assume that
$\a$ and $\a+\be$ appear consecutively in $\br = (r_1,\dots,r_n)$
(apply a sequence of short braid moves to $\br$ if necessary). In
other words, we may assume the existence of an integer $k$
satisfying $r_k= \a$ and $r_{k+1} = \a + \be$. Let $(s_1,\dots,
s_n)$ be the reduced expression corresponding to $\br$.  Parse
this reduced expression as $(\bu,s_{n-k},s_{n-k+1},\bv)$, and
denote the images of $\bu$ and $\bv$ in $W$ by $u$ and $v$,
respectively. For notational simplicity, write $s$ for $s_{n-k}$
and $t$ for $s_{n-k+1}$. Then $\a = v^{-1}(\a_t)$ and $\a + \be =
v^{-1} t (\a_s)$.  Note for future reference that $t(\a_s) =
s(\a_t) = \a_s+\a_t$.

Contractibility of $T$ will follow if we can produce a root
sequence for $ust$ of the form $ (v(\a), \ v(\a + \be), \ v(\be),
\ \ldots) = (\a_t, \ t(\a_s) = \a_s + \a_t, \ \a_s, \ \ldots)$.
Observe that
$$ v^{-1}(\a_s) = v^{-1}t(\a_s) - v^{-1}(\a_t) =
(\a+\be) - \a = \be \in \Phi(w) ,$$ hence $\a_s \in \Phi(u s t)$,
and consequently $\a_t = st(\a_s) \in \Phi(u)$. Thus, $u$ has a
reduced expression with last entry $t$, and this gives us a
reduced expression for $u s t$ of the form $(\dots , t, s, t)$.
The latter corresponds to a root sequence with the required
property.

The proof of (iii) $\Rightarrow$ (i) proceeds along similar lines.
We start with a root sequence $\br=(r_1,\dots,r_n)$ for $w$ for
which there is an index $k$ satisfying $r_k = \a+\be$ and
$r_{k+1}=\be$. We parse the reduced expression $(s_1,\dots,s_n)$
corresponding to $\br$ as $(\bu,s_{n-k},s_{n-k+1},\bv)$, denote
the images of $\bu$ and $\bv$ in $W$ as $u$ and $v$, and write $s$
for $s_{n-k}$ and $t$ for $s_{n-k+1}$. Then $\a+\be =
v^{-1}(\a_t)$ and $\be = v^{-1}t(\a_s)$.

We will be finished if we can show that $v$ has a reduced
expression with first entry $s$, say $(s,\bx)$. For then, letting
$x$ denote the image in $W$ of $\bx$, we would have
$$x^{-1}(\a_s) = -x^{-1}sts(\a_t) = -v^{-1}ts(\a_t) = v^{-1}(\a_t) -
v^{-1}t(\a_s)= (\a+\be) - \be = \a,$$ showing that $\a, \a + \be,
\be$ appear consecutively in some root sequence for $w$. An
equivalent requirement is for $v^{-1}(\a_s)$ to be a negative
root, and this follows because $$ v^{-1}(\a_s) = v^{-1}(\a_s+\a_t)
- v^{-1}(\a_t) = v^{-1}t(\a_s) - v^{-1}(\a_t) = \be - (\a+\be) =
-\a.$$ The proof is complete. \qed\enddemo

The next proposition will be needed in \S4.1.

\proclaim{Proposition 3.2.2} Let $w \in W$. Let $\a, \be \in
\Phi(w)$ be distinct and nonorthogonal. Suppose that there exist
root sequences $\br$ and $\br'$ for $w$ such that $\a \leq_{\br}
\be$ and $\be \leq_{\br'} \a$. Then there is a contractible
inversion triple of $w$ containing both $\a$ and $\be$.
\endproclaim

\demo{Proof} The relation $\a \leq_\br \be$ implies that in $\br$,
the root $\a$ appears to the left of $\be$; likewise, in $\br'$
the root $\be$ appears to the left of $\a$.  Now, by Remark 3.1.3,
$\br$ can be transformed into $\br'$ by a sequence of long and
short braid moves.  Since $\a$ and $\be$ are nonorthogonal, they
cannot both be involved in the same short braid move.  Hence,
there is a long braid move in our sequence of moves that involves
both $\a$ and $\be$. This implies the existence of a contractible
inversion triple of $w$ containing $\a$ and $\be$. \qed\enddemo

\head 4. The statistic $N(w)$ and free braidedness \endhead

Throughout this section, $(W,S)$ will be of arbitrary simply laced
type.

\subhead 4.1 The map $F_w$ \endsubhead

Let $w \in W$. Fix an arbitrary antisymmetric relation $\preceq$
on $\Phi(w)$ with the property that any two roots in $\Phi(w)$ are
comparable relative to $\preceq$.

Let ${\Cal C}(w)$ and ${\Cal I}(w)$ denote the set of commutation
classes of root sequences for $w$ and the set of contractible
inversion triples of $w$, respectively.  We define a map
$$ F_w : {\Cal C}(w) \longrightarrow
\{0,1\}^{{\Cal I}(w)},$$ depending on $\preceq$, as follows. Let
$C\in {\Cal C}(w)$ and let $\leq_C$ be the partial ordering of
$\Phi(w)$ determined by $C$; this is well-defined by Proposition
3.1.5.  Given any $\{\a,\a+\be,\be\}\in {\Cal I}(w)$, we define
$F_w(C)(\{\a, \a+ \be, \be \})$ to be $0$ if $\a$ and $\be$ are in
the same relative order with respect to $\leq_C$ and $\preceq$,
and otherwise we define $F_w(C)(\{\a, \a+ \be, \be\})$ to be $1$.
(Note that $\a$ and $\be$ are comparable under $\leq_C$.)

\proclaim{Theorem 4.1.1} The map $F_w$ is injective.
\endproclaim

\demo{Proof} Let $C$ and $C'$ be distinct root sequence
commutation classes of $w$, and let $\leq_C$ and $\leq_{C'}$ be
the respective partial orderings of $\Phi(w)$.  By Proposition
3.1.5, these partial orderings do not agree. Hence, there exist
nonorthogonal roots $\a, \be \in \Phi(w)$ such that $\a <_C \be$
and $\be <_{C'} \a$ (note that any two nonorthogonal roots in
$\Phi(w)$ are comparable under any given partial ordering arising
from a commutation class of $w$).  We now invoke Proposition 3.2.2
to deduce that $\a$ and $\be$ belong to some contractible
inversion triple $T$ of $w$. It is clear from the definition of
$F_w$ that $F_w(C)(T) \neq F_w(C')(T)$. \qed\enddemo

Theorem 4.1.1 has the following immediate

\proclaim{Corollary 4.1.2} Every $w\in W$ has at most $2^{N(w)}$
commutation classes, where $N(w)$ denotes the number of
contractible inversion triples of $w$. \qed\endproclaim

\remark{Remark 4.1.3} In \cite{4}, Elnitsky exhibits, in types
$A$, $B$ and $D$, a bijection from the set of all commutation
classes of an arbitrary group element to a set of rhombic tilings
of a polygon determined by that element.  B\'edard constructs in
\cite{1} a bijection from the set of all commutation classes of an
arbitrary Weyl group element $w$ to a certain set of functions
from $\Phi(w)$ to the set of positive integers.
\endremark

\subhead 4.2 Freely braided elements and commutation classes
\endsubhead

In this section, we prove that the freely braided elements of
Definition 2.2.1 achieve the bound of Corollary 4.1.2.

Given any $v\in V$, we denote by $v^\bot$ the subspace $\{w\in V :
B(v,w)=0\}$ of $V$.

\proclaim{Lemma 4.2.1} Let $w \in W$ have length at least $3$, and
let $\a_1, \a_2, \a_3 \in \Phi(w)$ be distinct. Suppose that $\br$
and $\br'$ are root sequences for $w$ such that $\a_2$ appears
between $\a_1$ and $\a_3$ in $\br$ but not in $\br'$. Then either
$\a_2^\bot \cap \{\a_1,\a_3\}\neq \emptyset$ or else there is a
contractible inversion triple of $w$ containing $\a_2$ and one of
$\a_1,\a_3$.
\endproclaim

\demo{Proof} By Remark 3.1.3, it is possible to transform $\br$
into $\br'$ using short and long braid moves.  Consideration of
such a sequence of moves yields the result. \qed\enddemo

\proclaim{Lemma 4.2.2} Let $w \in W$ be freely braided, and let $T
= \{\a,\a+\be,\be\}$ be a contractible inversion triple of $w$.
Suppose that $$
(\dots,\a,\g_1,\g_2,\dots,\g_k,\a+\be,\d_1,\d_2,\dots,\d_m,\be,\dots)
$$ is a root sequence for $w$.  Then each root $\g_i$ is orthogonal
to $\a$, and each root $\d_i$ is orthogonal to $\be$.
\endproclaim

\demo{Proof} First assume toward a contradiction that $\g_1$ is
not orthogonal to $\a$. We apply Lemma 4.2.1 with $\a_1 = \a$,
$\a_2 = \g_1$ and $\a_3 = \a + \be$. Because $T$ is contractible,
there is a root sequence for $w$ in which $\g_1$ does not appear
between the roots $\a$ and $\a + \be$. Further, since $w$ is
freely braided, and since $\g_1$ does not lie in $T$, it cannot be
the case that $\g_1$ belongs to a contractible inversion triple
containing either $\a$ or $\a+\be$. Lemma 4.2.1 then implies
$\g_1^\bot \cap \{\a,\a+\be\} \neq \emptyset$. But $\g_1$ is not
orthogonal to $\a$ by assumption; hence, $\g_1 \bot (\a + \be)$.
This in turn implies that $\g_1$ is not orthogonal to $\be$. Apply
Lemma 4.2.1 again with $\a_1 = \a$, $\a_2 = \g_1$ and $\a_3 =
\be$, and use the same reasoning to show $\g_1 \bot \be$.  This is
a contradiction.

An induction on $k$ now shows that the $\g_i$ are all orthogonal
to $\a$, and an analogous argument handles the $\d_i$.
\qed\enddemo

The following theorem and its corollary provide some justification
for the terminology ``freely braided element".

\proclaim{Theorem 4.2.3} Let $w \in W$ be freely braided. Then
there is a root sequence $\br$ for $w$ such that the roots in any
given contractible inversion triple of $w$ appear consecutively in
$\br$. Furthermore, every root sequence for $w$ is commutation
equivalent to such a root sequence $\br$.
\endproclaim

\demo{Proof} Let $\br'$ be an arbitrary root sequence for $w$, and
let $T = \{\a, \a + \be, \be\}$ be a contractible inversion triple
of $w$. By Remark 2.2.2 and symmetry, we may assume that in
$\br'$, the root $\a$ appears to the left of $\a + \be$, which in
turn appears to the left of $\be$.  According to Lemma 4.2.2 and
Proposition 3.1.1(a), there is a sequence of short braid moves
that transforms $\br'$ into a root sequence $\br''$ for $w$ in
which $\a$, $\a + \be$ and $\be$ appear consecutively in that
order.

The conclusion of the theorem follows from the above by an
induction on the number of contractible inversion triples. The
crucial point to notice is that if $T'$ is a consecutively
occurring contractible inversion triple in $\br'$, then the
sequence of short braid moves of the preceding paragraph may be
chosen so that the roots occurring in $T'$ remain consecutive in
$\br''$. This can be arranged because $w$ is freely braided.
\qed\enddemo

\proclaim{Corollary 4.2.4} If $w\in W$ is freely braided, then $w$
has precisely $2^{N(w)}$ commutation classes of root sequences,
where $N(w)$ is the number of contractible inversion triples of
$w$.
\endproclaim

\demo{Proof} In view of Corollary 4.1.2, it suffices to show that
$w$ has at least $2^{N(w)}$ commutation classes. Theorem 4.2.3
guarantees the existence of a root sequence $\br$ for $w$ in which
the elements of each contractible inversion triple occur
consecutively; moreover, the contractible inversion triples in
$\br$ are disjoint as $w$ is freely braided. Consequently, $F_w$
is surjective, and the conclusion follows. \qed\enddemo

The above corollary translates via Remark 3.1.3 into an equivalent
statement about commutation classes of reduced expressions for
$w$.

\subhead 4.3 Commutation graphs \endsubhead

Let $w \in W$. The {\it commutation graph} of $w$, $G(w)$, is
defined to be the graph with vertex set ${\Cal C}(w)$ and with
edge set consisting of all pairs $\{C,C'\}$ of commutation classes
with the following property: there exist representatives $\br \in
C$ and $\br' \in C'$ such that $\br$ and $\br'$ differ only by a
single long braid move.

We define the {\it parity} of each commutation class $C\in {\Cal
C}(w)$ to be the number $(-1)^{M(C)}$, where $M(C) = \sum_{T \in
{\Cal I}(w)} F_w(C)(T)$. Note that the parity of a commutation
class depends on the map $F_w$, which in turn depends on the
relation $\preceq$.

\proclaim{Proposition 4.3.1} Let $w \in W$. Then the graph $G(w)$
is bipartite.
\endproclaim

\demo{Proof}  Let $\{C,C'\}$ be an edge of $G(w)$.  By the
definition of edge and Remark 2.2.2, there is precisely one $T\in
{\Cal I}(w)$ such that $F_w(C)(T) \neq F_w(C')(T)$. Hence, the
classes $C$ and $C'$ have unequal parity. \qed\enddemo

\remark{Remark 4.3.2} Let $\Gamma$ be the Coxeter graph of
$(W,S)$. If $\Gamma$ itself is bipartite, then Proposition 4.3.1
can be established using a simple argument counting the number of
occurrences of each of the two types of generators in a reduced
expression for $w$.
\endremark

\head 5. Type $A$ \endhead

For the remainder of this paper, we shall assume that $(W,S)$ is
of type $A$.  It will be convenient to work with a particular
realization of $(W,S)$, under which $W$ is the group ${\Cal S}_n$
of permutations of $\{1,\dots,n\}$, and $S$ is the set of simple
transpositions $s_1=(1,2),\dots,s_{n-1}=(n-1,n)$. Let $\e_1,\dots
,\e_n$ be the standard basis of ${\Bbb R}^n$, and let $\langle
\cdot,\cdot \rangle$ be the inner product on ${\Bbb R}^n$
satisfying $\langle \e_i,\e_j \rangle = \d_{ij}$.  Allow ${\Cal
S}_n$ to act on ${\Bbb R}^n$ by permuting the indices of
$\e_1,\dots,\e_n$, and set $\a_i = \e_i-\e_{i+1}$
$(i=1,\dots,n-1)$. Then the vectors $\a_1,\dots,\a_{n-1}$ may be
regarded as the simple roots corresponding to $s_1,\dots,s_{n-1}$,
and the form $\frac{1}{2}\langle \cdot,\cdot \rangle$ may be
regarded as the Coxeter form $B$.  With these identifications, we
have $\Phi = \{ \e_i-\e_j : i\neq j\}$ and $\Phi^+ =
\{\e_i-\e_j:i<j\}$.

Let $w\in {\Cal S}_n$.  By the ``1-line notation" for $w$, we
shall mean the string of numbers $w(1)\cdots w(n)$.  Let $(i,j)$
be a (not necessarily simple) transposition in ${\Cal S}_n$, with
$i<j$. Note that the 1-line notation for the product $w(i,j)$ can
be obtained from that for $w$ by interchanging the numbers in
positions $i$ and $j$. Further, we have $\ell(w(i,j))<\ell(w)
\Leftrightarrow w(i)>w(j) \Leftrightarrow \e_i-\e_j \in \Phi(w)$.

The set of inversion triples of any permutation $w\in {\Cal S}_n$
can be read off of the 1-line notation for $w$. Specifically, a
triple $\{\e_i-\e_j,\e_i-\e_k,\e_j-\e_k\}$ of positive roots is an
inversion triple of $w$ if and only if $w(i)>w(j)>w(k)$.

\subhead 5.1 Contractibility in type $A$ \endsubhead

The following proposition stands in contrast to Example 2.2.3.

\proclaim{Proposition 5.1.1} Let $w\in {\Cal S}_n$.  Then every
inversion triple of $w$ is contractible.
\endproclaim

\demo{Proof} We argue by induction on $m=\ell(w)$.  The assertion
is vacuously true if $0\leq m<3$. Suppose that $m\geq 3$ and that
$w$ has at least one inversion triple. Let
$T=\{\e_i-\e_j,\e_i-\e_k,\e_j-\e_k\}$ be an inversion triple of
$w$.

Consider the situation where there is a positive integer $l<n$
such that $\{l,l+1\}\cap \{i,j,k\}=\emptyset$ and $w(l)>w(l+1)$.
Then $\ell(ws_l)<\ell(w)$, and $T$ is an inversion triple of
$ws_l$. By the inductive hypothesis, $T$ is contractible relative
to $ws_l$.  If $(r_1,\dots,r_{m-1})$ is any root sequence for
$ws_l$ in which the elements of $T$ appear consecutively, then
$(\a_l,s_l(r_1),\dots,s_l(r_{m-1}))$ is a root sequence for $w$ in
which the elements of $T$ appear consecutively, since $s_l$ fixes
the elements of $T$ pointwise.

We may therefore assume that $w(l)<w(l+1)$ for every positive
integer $l<n$ such that $\{l,l+1\}\cap \{i,j,k\}=\emptyset$.

Suppose that $i<j-1$.  Since $w(i)>w(j)$, the previous paragraph
implies that either $w(i)>w(i+1)$ or $w(j-1)>w(j)$.  Suppose for
the moment that $w(i)>w(i+1)$. Then $\ell(ws_i)<\ell(w)$, and the
set $T'=\{\e_{i+1}-\e_j,\e_{i+1}-\e_k,\e_j-\e_k\}$ is an inversion
triple of $ws_i$.  By the inductive hypothesis, $T'$ is
contractible relative to $ws_i$.  If $(r_1,\dots,r_{m-1})$ is any
root sequence for $ws_i$ in which the elements of $T'$ appear
consecutively, then $(\a_i,s_i(r_1),\dots,s_i(r_{m-1}))$ is a root
sequence for $w$ in which the elements of $T$ appear
consecutively.

Similar reasoning handles the case where $w(j-1)>w(j)$. The
induction step is therefore proved for $i<j-1$.

If $j<k-1$, then one argues in a similar way. It remains to deal
with the case where $j=i+1$ and $k=j+1$. Here, the element $w$ has
a reduced expression of the form $(\dots,s_i,s_{i+1},s_i)$, and it
is clear that the first three entries of the corresponding root
sequence belong to $T$. \qed\enddemo

Let $w\in {\Cal S}_n$ and let $v\in {\Cal S}_k$, where $k\leq n$.
Suppose that there do not exist $i_1,\dots,i_k\in \{1,\dots,n\}$
with $i_1<\cdots <i_k$ such that the numbers $w(i_1),\dots,w(i_k)$
are in the same relative order as $v(1),\dots,v(k)$.  Then one
says that $w$ ``avoids the pattern" $v(1)\cdots v(k)$.  It can be
deduced from Proposition 5.1.1 and the discussion that precedes it
that a given permutation is freely braided if and only if it
avoids each of the four patterns $3421, 4231, 4312$ and $4321$.

\subhead 5.2 Converse to Corollary 4.2.4 in type $A$ \endsubhead

Let $w\in {\Cal S}_n$.  Let ${\Cal T}$ be a nonempty subset of
${\Cal I}(w)$, the set of all (contractible) inversion triples of
$w$. We say that the map $F_w$ of \S4.1 {\it separates} ${\Cal T}$
if every map from ${\Cal T}$ to $\{0,1\}$ is the restriction of
some element of $F_w({\Cal C}(w))$. Clearly, if $F_w$ fails to
separate a nonempty subset of ${\Cal I}(w)$, then $F_w$ is not
surjective.

The following theorem is a partial converse to Corollary 4.2.4.
Its proof relies on Proposition 5.1.1, which, as was observed in
Example 2.2.3, does not hold for all simply laced Coxeter systems.

\proclaim{Theorem 5.2.1} Suppose that $w\in {\Cal S}_n$ has
$2^{N(w)}$ inversion triples.  Then $w$ is freely braided.
\endproclaim

\demo{Proof} Assume toward a contradiction that there is a
non-freely-braided $w\in {\Cal S}_n$ with $2^{N(w)}$ inversion
triples (all of which must be contractible by Proposition 5.1.1).
In other words, $w$ is non-freely-braided and $F_w$ is surjective.
Let $\e_i-\e_j$ be an element of $\Phi(w)$ that belongs to at
least two inversion triples of $w$, and assume that $j-i$ is as
large as possible relative to this property.  Let $T$ and $T'$ be
distinct inversion triples of $w$ containing $\e_i-\e_j$. By
symmetry, there are six cases to consider.

Case 1. $\e_j-\e_k\in T$ and $\e_j-\e_l\in T'$, where $i<j<k<l$.

\noindent Observe that $\e_k-\e_l \notin \Phi(w)$. Otherwise, the
root $\e_i-\e_l$ would belong to at least two inversion triples of
$w$, namely $\{\e_i-\e_k,\e_i-\e_l,\e_k-\e_l\}$ and
$\{\e_i-\e_j,\e_i-\e_l,\e_j-\e_l\}$, contradicting our choice of
$\e_i-\e_j$.

By Proposition 5.1.1 and Remark 2.2.2, there is a root sequence
for $w$ of the form
$$(\dots,\e_i-\e_j,\e_i-\e_k,\e_j-\e_k,\dots,\e_i-\e_l,\dots,
\e_j-\e_l,\dots)\,\,\text{ or}$$
$$(\dots,\e_j-\e_l,\dots,\e_i-\e_l,\dots,\e_i-\e_j,\e_i-\e_k,
\e_j-\e_k,\dots).$$ Since $\e_k-\e_l \notin \Phi(w)$, and since
$\e_i-\e_k$ and $\e_i-\e_l$ are not orthogonal, these last two
roots maintain their positions relative to one another in every
root sequence for $w$. If $\e_i-\e_k$ lies always to the left of
$\e_i-\e_l$, then $\e_i-\e_j$ cannot at the same time lie to the
left of $\e_i-\e_k$ and to the right of $\e_i-\e_l$.  A similar
statement can be made if $\e_i-\e_k$ lies always to the right of
$\e_i-\e_l$. It follows that $F_w$ does not separate $\{T,T'\}$.

Case 2. $\e_k-\e_i\in T$ and $\e_l-\e_i\in T'$, where $k<l<i<j$.

\noindent One argues as in Case 1.

Case 3. $\e_k-\e_i\in T$ and $\e_j-\e_l\in T'$, where $k<i<j<l$.

\noindent Note that $\e_k-\e_j \in T\subseteq \Phi(w)$. Since both
$\e_k-\e_i$ and $\e_j-\e_l$ also lie in $\Phi(w)$, the root
$\e_k-\e_j$ belongs to at least two inversion triples of $w$,
namely $\{\e_k-\e_j,\e_k-\e_l,\e_j-\e_l\}$ and
$\{\e_k-\e_i,\e_k-\e_j,\e_i-\e_j\}$.  This contradicts our choice
of $\e_i-\e_j$.

Case 4. $\e_i-\e_k \in T$ and $\e_j-\e_l\in T'$, where $i<k<j<l$.

\noindent We have $\e_k-\e_j\in T\subseteq \Phi(w)$. Hence,
$\e_k-\e_l = (\e_k-\e_j)+(\e_j-\e_l)\in \Phi(w)$. It follows that
$\e_i-\e_l$ belongs to at least two inversion triples of $w$,
namely $T'$ and $\{\e_i-\e_k,\e_i-\e_l,\e_k-\e_l\}$, contradicting
our choice of $\e_i-\e_j$.

Case 5. $\e_k-\e_i\in T$ and $\e_l-\e_j\in T'$, where $k<i<l<j$.

\noindent One argues as in Case 4.

Case 6.  $\e_i-\e_k\in T$ and $\e_l-\e_j\in T'$, where $i<k<l<j$.

\noindent Here, we have $T=\{\e_i-\e_k,\e_i-\e_j,\e_k-\e_j\}$ and
$T'=\{\e_i-\e_l,\e_i-\e_j,\e_l-\e_j\}$. Suppose first that
$\e_k-\e_l\notin \Phi(w)$. Then the roots $\e_i-\e_k$ and
$\e_i-\e_l$, which are not orthogonal, maintain their positions
relative to one other in every root sequence for $w$.  It follows
that $F_w$ does not separate $\{T,T'\}$.

Finally, suppose that $\e_k-\e_l \in \Phi(w)$.  Then the sets $T''
=\{\e_k-\e_l,\e_k-\e_j,\e_l-\e_j\}$ and $T''' = \{\e_i-\e_k,
\e_i-\e_l, \e_k-\e_l\}$ are both inversion triples of $w$.  We
claim that $F_w$ does not separate $\{T,T'',T'''\}$.  To see this,
observe that if $C$ is a commutation class of $w$ relative to
which $\e_i-\e_k$ lies to the left of $\e_k-\e_j$ (this determines
$F_w(C)(T)$) and $\e_k-\e_j$ lies to the left of $\e_k-\e_l$ (this
determines $F_w(C)(T'')$), then $\e_i-\e_k$ lies to the left of
$\e_k-\e_l$ (so that $F_w(C)(T''')$ is also determined).
\qed\enddemo

\noindent {\bf Note added in proof.} The authors have found a
proof that the converse of Corollary 4.2.4 holds for arbitrary
simply laced Coxeter systems.

\leftheadtext{} \rightheadtext{}

\end